\documentclass{amsart}

\usepackage{amssymb}

\usepackage{graphicx}

\usepackage{amscd}
\usepackage[cmtip,all]{xy}

\newtheorem{theorem}{Theorem}[section]

\newtheorem{proposition}[theorem]{Proposition}
\theoremstyle{definition}

\newtheorem{corollary}[theorem]{Corollary}

\theoremstyle{remark}

\numberwithin{equation}{section}

\begin{document}

\title[$K(s)^*(B\mathcal{G})$ for some Frobenius complements]{ Mod 2 Morava $K$-theory for Frobenius complements of exponent dividing $2^n \cdot 9$ }


\author{Malkhaz Bakuradze}
\address{ Tbilisi State university, Georgia}
\curraddr{}
\email{malkhaz.bakuradze@tsu.ge}
\thanks{Research was supported by Volkswagen Foundation, Ref.: I/84 328 and by GNSF, ST08/3-387}


\subjclass[2000]{55N20; 55R12; 55R40}

\date{}

\dedicatory{}

\commby{}


\begin{abstract}

 We determine the cohomology rings $K(s)^*(B\mathcal{G})$ at 2 for all finite Frobenius complements $\mathcal{G}$ of exponent dividing $2^n \cdot 9$.
\end{abstract}
\maketitle



Let $V$ be an abelian group, and let $\mathcal{G}$ be a group of automorphisms of $V$. If $\mathcal{G}$ has exponent $2^n \cdot 3^k $ for
$0 \leq n $ and $0 \leq k \leq 2 $ and $\mathcal{G}$ acts freely on $V$, then $\mathcal{G}$ is finite (see \cite{JM} Theorem 1.1). Every finite group that acts freely on an abelian group is isomorphic to a Frobenius complement in some finite Frobenius group (see \cite{JM} Lemma 2.6).
By the classification of finite Frobenius complements (see \cite{P}) the quotient of $\mathcal{G}$ by its maximal normal 3-subgroup $\mathcal{H}$ is isomorphic to a cyclic 2-group $\mathcal{C}$, a generalized quaternion group $Q$, the binary tetrahedral group $2\mathcal{T}$ of order 24 (or SL(2,3)),  or the binary octahedral group $2\mathcal{O}$ of order 48. Then Atiyah-Hirzebruch-Serre spectral sequence for $\mathcal{H} \lhd \mathcal{G}$ implies
 that at 2 the ring $K(s)^*(B\mathcal{G})$ is isomorphic to $K(s)^*(B\mathcal{K})$, for $\mathcal{K}=\mathcal{G}/\mathcal{H}$ is either $\mathcal{C}, Q, 2\mathcal{T}, 2\mathcal{O}.$
 For the cyclic group $\mathcal{C}=\mathbb{Z}/2^{k}$, $K(s)^*(B\mathbb{Z}/2^{k})=\mathbb{F}_2[v_s,v_s^{-1}][u]/(u^{2^{ks}}).$ For the generalized quaternion group $Q_{2^{m+2}}$ we have Theorem 1.1 of \cite{BV}. We deduce Morava $K$-theory rings at 2 for the groups $2\mathcal{T}$ and $2\mathcal{O}$ as certain subgroups in $K(s)^*(BQ_8)$ and $K(s)^*(BQ_{16})$ respectively (Proposition \ref{prop:T} and Proposition \ref{prop:2O}.)

In \cite{BP} we proved the following formula for the first Chern class of the transferred line complex bundle:  Let $X \rightarrow Y$ be the regular two covering defined by free action of $\mathbb{Z}/2$ on $X$ and let $\theta \rightarrow Y$ be the associated line complex bundle;
Let $\xi \rightarrow X$ be a complex line bundle and let $\zeta \rightarrow Y$ be the plane
bundle, transferred from
$\xi$ by Atiyah transfer \cite {At}.
Then for $Tr^*: K(s)(X)\rightarrow K(s)^*(Y)$, the transfer homomorphism \cite{Ada} for our covering $X \rightarrow Y$, one has

\begin{equation}
\label{eq:tr}
Tr^*(c_1(\xi))=c_1(\theta)+c_1(\zeta)+v_s\sum_{i=1}^{s-1}c_1(\theta)^{2^s-2^{i}}c_2(\zeta)^{2^{i-1}}.
\end{equation}

\bigskip
We show that formula \ref{eq:tr} plays major role in the ring structure $K(s)^*(B\mathcal{G})$ at 2 for aforementioned groups and gives another derivations for some related rank one Lie groups.
\bigskip

Much of our note is written in terms of Theorem 1.1 of \cite{BV}. Let
$$
G= \langle a,b \ | a^{2^{m+1}}=1, \ b^2=a^e, bab^{-1}=a^r \rangle,
\ \ m \geq 1
$$
and either $e=0,$ $r=-1$ (the dihedral group $D_{2^{m+2}}$ of
order $2^{m+2}$), $e=2^m,$ $r=-1$ (the generalized quaternion
group $Q_{2^{m+2}}$) or $m\geq 2,$ $e=0,$ $r=2^m-1$ (the
semidihedral group $SD_{2^{m+2}}$).

Spectral sequence consideration (see \cite{TY}) imply that $K(s)(BG)$ is generated by following Chern classes
$|c|=|x|=2,$ $|c_2|=4$:

\begin{align*}
c&=c_1(\eta_1), \ \eta_1:G/\langle a \rangle \cong\mathbf
{Z}/2\rightarrow \mathbb{C}^*, \ b\mapsto -1; \\
x&=c_1(\eta_2), \ \eta_2:G/\langle a^2,b \rangle
\cong\mathbf{Z}/2\rightarrow\mathbb {C}^*, \ a\mapsto -1;
\end{align*}

and $c_2=c_2(\xi_{\pi_{1}})$, where
$\xi_{\pi_{1}}\rightarrow B\langle a,b \rangle$
is the plane bundle transferred from the canonical line bundle
$\xi \rightarrow B\langle a \rangle ,$
for the double covering
$
\pi_{1}:B\langle a \rangle \rightarrow B\langle a,b \rangle
$
corresponding to $\eta_1$.

The ring structure is the result of the formula for transferred
first Chern class \ref{eq:tr}. See \cite{BV}.

\bigskip

Let $N$ be the normalizer of $U(1)$ in $S^3$.
The normalizes of the maximal torus in $SO(3)$ is $O(2)=U(1)\rtimes \mathbb{Z}/2$ and
$\mathbb{Z}/2$ acts on $K(s)^*BU(1)=K(s)^*[[u]]$ by  $[-1]_{F}(u)$ as above.

Since $BU(1)\hat{}_p=[colim_n B \mathbb{Z}/(p^n)]\hat{}_p$, we have
$K(s)^*(BO(2))=K(s)^*(lim_m(BD_{2^{m+2}}))=K(s)^*(lim_m(BSD_{2^{m+2}}))$ and $K(s)^*(BN)=K(s)^*(lim_m(BQ_{2^{m+2}}))$.

\medskip

Thus Theorem 1.1 of \cite{BV} implies

\begin{corollary}
\label{prop:BO(2)}
$K(s)^*(BO(2))=K(s)^*[[c,c_2]]/(c^{2^s}, v_s c
\sum_{i=1}^{s}c^{2^s-2^i}c_2^{2^{i-1}}),$ where  $c=c_1(det \eta)$ and
$c_2=c_2(\eta)$ are the Chern classes of the bundle $\eta \rightarrow BO(2)$, the complexification of canonical $O(2)$ bundle.
\end{corollary}

\bigskip

\begin{corollary}
\label{prop:CPm}
$K(s)^*(BN)=K^*(s)[[c,c_2]]/(c^{2^s}, c^2+v_s c
\sum_{i=1}^{s}c^{2^s-2^i}c_2^{2^{i-1}}),$

where $c=c_1(\nu)$ is the Chern class  of $\nu$ the pullback bundle of the canonical real line bundle by
$N \rightarrow N/U(1)=\mathbb{Z}/2$ and  $c_2=c_2(p^*(\zeta))$ is the Euler class of the pullback bundle of the canonical quaternionic line bundle by the inclusion $N \subset S^3.$
\end{corollary}

Then $RP^2\rightarrow BO(2)\rightarrow BO(3)$ is the projective bundle of the canonical $SO(3)$ bundle. Hence the pullback of the complexification  of this canonical $SO(3)$ bundle splits over $BO(2)$ as $\eta \oplus det\eta$. Note that $c_1(det \eta)=c_1(\eta)+v_sc_2(\eta)^{2^{s-1}}$ modulo transfer for the covering $BU(1)\rightarrow BO(2).$
Thus $K(s)^*(BSO(3))$ is subring in $K(s)^*(BO(2))$ generated by $v=c^2+v_scc_2^{2^{s-1}}+c_2$ and $w=cc_2.$
This implies

\begin{corollary}
\label{cor:BSO(3)}
$K(s)^*(BSO(3))=K(s)^*[[v,w]](f_s(v,w), g_s(v,w)),$

where $|v|=4$, $|w|=6$, and $f_s=f_s(v,w)$, $g_s=g_s(v,w)$ are determined by $f_2=vw$, $g_2=w^2$ and for $s>2$
\end{corollary}

\medskip

$
f_s=\begin{cases}
f^2_{s-1}& \text{$s$ even},\\
\frac{f_{s-1}g_{s-1}}{v}+wv^{2^{s-1}-1} & \text{$s$ odd},
\end{cases}
$

\medskip

$
g_s =\begin{cases}
g^2_{s-1}& \text{$s$ odd},\\
\frac{f_{s-1}g_{s-1}}{v}+wv^{2^{s-1}-1} & \text{$s$ even}.
\end{cases}
$

\bigskip

Our main result is the following.

Let $\mathcal{G}$ be a group acting freely on an abelian group. Let $\mathcal{G}$ be of exponent dividing $2^n\cdot 9$ (hence $\mathcal{G}$ is necessarily finite, as above) and let $\mathcal{H}\lhd \mathcal{G}$ be the maximal normal 3-subgroup.

\begin{theorem}
\label{theorem2}
As a ring $K(s)^*(B\mathcal{G})$ has one of the following forms

\medskip

(i) If $\mathcal{G}/\mathcal{H}$=$Q_8$ , then $K(s)^*(B\mathcal{G})=K(s)^*[c,x,c_2]/R$

and the relations $R$ are determined by

\medskip

$
c^{2^s}=x^{2^s}=0,
$
$
v_scc_2^{2^{s-1}}=v_s \sum_{i=1}^{s-1}c^{2^s-2^i+1}c_2^{2^{i-1}}+c^2,
$
$
v_s^2c_2^{2^s}=c^2+cx+x^2,
$
$
v_sxc_2^{2^{s-1}}=v_s\sum_{i=1}^{s-1}x^{2^s-2^i+1}c_2^{2^{i-1}}+x^2.
$

\bigskip

(ii) If $\mathcal{G}/\mathcal{H}$=$Q_{2^{m+2}}$ , $m>1$ , then $K(s)^*(B\mathcal{G})=K(s)^*[c,x,c_2]/R,$

and the relations $R$ are determined by

\medskip

$
c^{2^s}=x^{2^s}=0,
$
$
v_scc_2^{2^{s-1}}=v_s \sum_{i=1}^{s-1}c^{2^s-2^i+1}c_2^{2^{i-1}}+c^2,
$
$
v_s^{2\kappa(m)}c_2^{2^{ms}}=cx+x^2,
$
$
v_sxc_2^{2^{s-1}}=v_sx\sum_{i=1}^{s-1}c^{2^s-2^i}c_2^{2^{i-1}}
+\sum_{i=1}^{ms}v_s^{1+\kappa(m)+2^{ms}-2^i}c_2^{(2^{ms}+1)2^{s-1}-(2^s-1)2^{i-1}}\\+cx,
$

where $\kappa(m)=\frac{2^{ms}-1}{2^s-1}.$

\bigskip

(iii) If $\mathcal{G}/\mathcal{H}$=$2\mathcal{T}$, then $K(s)^*(B\mathcal{G})=K(s)^*[c_2]/c_2^{(2^{s}+1)2^{s-1}}.$

\bigskip

(iv) If $\mathcal{G}/\mathcal{H}$=$2\mathcal{O}$, then

$K(s)^*(B\mathcal{G})=K(s)^*[c, c_2]/ (c^{2^s},
c^2+v_sc\sum_{i=1}^{s}c^{2^s-2^i}c_2^{2^{i-1}},
c_2^{(2^{s}+1)2^{s-1}}).$

\bigskip

(v) If $\mathcal{G}/\mathcal{H}$=$\mathbb{Z}/2^k$, then $K(s)^*(B\mathcal{G})=K(s)^*[c]/c^{2^{ks}}.$

\medskip

Here in all cases $|c|=|x|=2$, $|c_2|=4.$

\end{theorem}

\bigskip

The statement (v) is clear. (i) and (ii) follow from Theorem 1.1 of \cite{BV} for $Q_8$ and $Q_{2^{m+2}}$ respectively.
What remains is to consider the cases of binary tetrahedral and binary octahedral groups.

\section{Binary Polyhedral groups}

As it is known any finite subgroup of $SO(3)$ is either a cyclic
group, a dihedral group or one of the groups of a Platonic solid:
tetrahedral group $\mathcal{T} \cong A_4$, cube/octahedral group
$\mathcal{O} \cong S_4,$ or icosahedral group $\mathcal{I} \cong
A_5 $. We consider the preimages of the latter groups under
the covering homomorphism $S^3\rightarrow SO(3)$.

\medskip

\subsection*{Binary tetrahedral group}

Binary tetrahedral group $2\mathcal{T}$ as the group of 24 units
in the ring of Hurwitz integers $2\mathcal{T}$ is given by
$
\{\pm 1, \pm i, \pm j, \pm k, \frac{1}{2}(\pm 1 \pm i \pm j \pm k
)\}.
$

This group can be written as a semidirect product $2T=Q_8 \rtimes
\mathbb{Z}/3$, where $Q_8$ is the quaternion group consisting of
the 8 Lipschitz units $\pm 1, \pm i, \pm j, \pm k $ and
$\mathbb{Z}/3$ is the cyclic group generated by $ -\frac{1}{2}(1+i+j+k ).$ The cyclic group acts on the normal
subgroup $Q_8$ by conjugation. So that the generator of
$\mathbb{Z}/3$ cyclically rotates $i,j,k.$

Consider now Morava $K$-theory at 2. Then relations of Theorem 1.1 of \cite{BV} for $K(s)^*(BQ_8)$ imply that its subring of
invariants under $\mathbb{Z}/3$ action is generated by $c_2$: the
generator of $\mathbb{Z}/3$ cyclically rotates $c$, $x$ and
$c+x+v_sc^{2^{s-1}}x^{2^{s-1}}$. If ignoring the powers of $v_s$
then the first and second elementary symmetric functions in these
three symbols are equal to $c_2^{2^{s-1}}$ and $c_2^{2^s}$
respectively and the third is zero. It follows that $K(s)^*(B2\mathcal{T})\cong
[K(s)^*(BQ_8)]^{\mathbb{Z}/3}.$

\begin{proposition}
\label{prop:T} $K(s)^*(B2\mathcal{T})\cong
K(s)^*[c_2]/c_2^{(2^{s}+1)2^{s-1}},$ where $|c_2|=4.$
\end{proposition}

\bigskip

\subsection*{Binary octahedral group $2\mathcal{O}$} This group is given
as the union of the 24 Hurwitz units
$\{\pm 1, \pm i, \pm j, \pm k, \frac{1}{2}(\pm 1 \pm i \pm j \pm k)\}$
with all 24 quaternions obtained from  $\frac{1}{\sqrt{2}}(\pm 1
\pm i+ 0j + 0k )$ by permutation of coordinates.

The generalized quaternion group $Q_{16}$ forms a subgroup of
$2\mathcal{O}$ and its conjugacy classes has 3 members. Therefore
by the transfer argument $B2\mathcal{O}$ is a stable wedge summand
of $BQ_{16}$ after localized at 2, meaning $K(s)^*(B2\mathcal{O})$
is the subring in $K(s)^*(BQ_{16})$ at 2. We show that
this is the subring generated by two symbols $c$ and $c_2$ of
Theorem 1.1 of \cite{BV}. Namely one has

\begin{proposition}
\label{prop:2O} $K(s)^*(B2\mathcal{O})$ is isomorphic to
$$K(s)^*[c, c_2]/ (c^{2^s},
c^2+v_sc\sum_{i=1}^{s}c^{2^s-2^i}c_2^{2^{i-1}},
c_2^{(2^{s}+1)2^{s-1}}),$$  where $|c|=2,$ $|c_2|=4.$
\end{proposition}

\subsection*{Binary icosahedral group} $2\mathcal{I}$ is given as the union of
the 24 Hutwitz units
$\{\pm 1, \pm i, \pm j, \pm k, \frac{1}{2}(\pm 1 \pm i \pm j \pm k)\}$
with all 96 quaternions obtained from  $\frac{1}{2}(0\pm 1\pm
i\pm\varphi^{-1}j\pm \varphi k )$ by even permutation of
coordinates. Here $\varphi=\frac{1}{2}(1+\sqrt{5})$ is the golden
ratio. This group is isomorphic to $SL_2(5)$-the group of all $2\times 2$
matrices over $\mathbb{F}_5$ with unit determinant.

Among other subgroups the relevant subgroup is the binary
tetrahedral group formed by Hurwitz units. Then coset
$2\mathcal{I}/ 2\mathcal{O}$ has 5 members hence by the transfer
argument again $B2\mathcal{I}$ splits off $B2\mathcal{O}$ after
localized at 2. Thus we obtain
\begin{equation*}
K(s)^*B(2\mathcal{I})\cong
K(s)^*B(2\mathcal{T}).
\end{equation*}
\bigskip

\bibliographystyle{amsplain}

\end{document}